\documentclass{amsproc}
\newcommand{\cal}{\mathcal }

\usepackage{graphicx, overpic}
\usepackage{natbib}
\usepackage{latexsym,array,delarray,amsmath,amsfonts,amssymb,amsthm}
\usepackage{bm}
\bmdefine{\Bt}{t}
\bmdefine{\BX}{X}
\bmdefine{\BY}{Y}
\bmdefine{\BZ}{Z}
\bmdefine{\BB}{B}
\bmdefine{\BM}{M}
\bmdefine{\BD}{D}
\bmdefine{\Bi}{i}
\bmdefine{\Bj}{j}
\bmdefine{\Bx}{x}
\bmdefine{\By}{y}
\bmdefine{\Bz}{z}
\bmdefine{\Bw}{w}
\bmdefine{\Ba}{a}
\bmdefine{\Bb}{b}
\bmdefine{\Bc}{c}
\bmdefine{\Bh}{h}
\bmdefine{\Bg}{g}
\bmdefine{\Bu}{u}

\def\sg{{\rm semigroup}}
\def\cone{{\rm cone}}

\def\Qsat{Q_{\rm sat}}

\newtheorem{thm}{Theorem}[section]
\newtheorem{lemma}[thm]{Lemma}
\newtheorem{prop}[thm]{Proposition}

\newtheorem{defn}[thm]{Definition}
\newtheorem{ex}[thm]{Example}

\newcommand{\Z}{{\mathbb Z}}
\newcommand{\N}{{\mathbb N}}

\newcommand{\R}{\mathbb R}

\def\comment#1{\textit{[#1]}}
\def\comment#1{}

\newcommand{\zz}{\mathbb{Z}}
\newcommand{\nn}{\mathbb{N}}

\newcommand{\rr}{\mathbb{R}}

\newcommand{\relint}{{\rm relint}}

\begin{document}
\title{
Saturation points on faces of a rational polyhedral cone}
\author{ Akimichi Takemura }
\address{Department of Mathematical Informatics,  
   University of Tokyo, Bunkyo, Tokyo, Japan.  
   takemura@stat.t.u-tokyo.ac.jp} 
\author{Ruriko  Yoshida} 
\address{ Department of Statistics, University of Kentucky, Lexington,
  KY USA. ruriko@ms.uky.edu}

\maketitle

\begin{abstract}
  Different commutative semigroups may have a common saturation.  We
  consider distinguishing semigroups with a common
  saturation based on their ``sparsity''.  We propose to qualitatively
  describe sparsity of a semigroup by considering which faces of the
  corresponding rational polyhedral cone have saturation points.  For
  a commutative semigroup we give a necessary and sufficient condition
  for determining which faces have saturation points.  We also show
  that we can construct a
  commutative semigroup with arbitrary consistent patterns of faces
  with saturations points.
\end{abstract}

\begin{quotation}
\par
\end{quotation}

\section{Introduction}
Let $A=\{\Ba_1,\ldots,\Ba_n\}$, $\Ba_i \in \Z^d$, $i=1,\dots,n$, be a
finite set of integral points and let $Q=Q(A)$ denote the
commutative semigroup generated by $\Ba_1,\ldots,\Ba_n$.  In our
previous paper (\cite{takemura-yoshida2006}) we studied properties of
{\em holes}, which are the difference between the semigroup and its
{\em saturation}. We 
gave some necessary and sufficient conditions for the finiteness of the set 
of holes. % (Theorem \ref{thm:previous-paper} below).  
In this paper we give a more detailed description of 
how holes are located, when there are infinitely many holes.  
Studying holes of $Q$
finds applications in many areas, such as commutative algebra
(\cite{BGT, BG, BGcover, KantorSarkaria, GKZ,sturmfels1996}),
optimization (\cite{aardaletal3,aardaletal2,aardaletal1}), number
theory (\cite{Barvinok2002}), and statistics
(\cite{dobra-karr-sanil2003, Cox00, Cox02,
  ohsugi-hibi-contingency-tables-2006, Vlach1986}).

Let $K=\cone(\Ba_1,\allowbreak\ldots,\Ba_n)$ be the rational
polyhedral cone generated by $\Ba_1,\ldots,\Ba_n$.  In this paper,
without essential loss of generality, we assume that the lattice
generated by $\Ba_1,\ldots,\Ba_n$ is $\Z^d$. In this case the 
saturation $\Qsat$ of $Q$ is defined by $\Qsat=K \cap \Z^d$.  The
elements of $H=\Qsat \setminus Q$ are called {\em holes} of $Q$.  We assume
that $K$ is a {\em pointed} cone with non-empty interior.
% throughout this paper.
Let $B$ denote the unique minimal {\em Hilbert basis} of $K$ (i.e.\ the unique
minimal generator of $K \cap \Z^d$).  In the following we simply say
the Hilbert basis instead of the unique minimal Hilbert basis.

All holes have to be close to the boundary of $K$.  However
holes may be concentrated only around some low dimensional proper
faces of $K$ or they may be distributed all over the boundary of $K$.
In the latter case $Q$ can be considered as more sparse than the
former case.  In order to qualitatively describe this distinction we
define the notion of an {\em almost saturated face} of $K$ and show how to
determine whether a given face $F$ is almost saturated or not.  In
this paper we assume that $K$ is a given rational polyhedral cone
and the purpose of our investigation is to differentiate semigroups
with common saturation $K \cap \Z^d$ by their configurations of the almost
saturated faces.

For the rest of this section we summarize some relevant
definitions and results from \cite{takemura-yoshida2006}.
% For a summary of definitions and notation from
% \cite{takemura-yoshida2006}, we have the followings: Let $Q$ be the
% semigroup generated by $\Ba_1,\ldots,\Ba_n$, let
%  and let $L$ be the lattice generated by
% $\Ba_1,\ldots,\Ba_n$.  Then the semigroup $\Qsat = K \cap L$ is called
% the {\em saturation} of the semigroup $Q$ and clearly $Q \subset
% \Qsat$.  We call $Q$ {\em saturated} if $Q = \Qsat$ (also this is
% called {\em normal}).  Let $H=\Qsat \setminus Q$ be the set of holes.
We call $\Ba \in \Qsat$, $\Ba\neq 0$, a {\em fundamental hole} if
$\Qsat \cap (\Ba + (-Q)) = \{ \Ba \}$.  Let $H_0$ be the set of all
fundamental holes in $Q$.  $H_0$ is always finite for any
given semigroup by Proposition 3.1 in \cite{takemura-yoshida2006}.
$\Ba \in Q$ is called a {\em saturation point} if 
$\Ba + \Qsat \subset Q$.  Let $S$ be the set of all saturation points of
the semigroup $Q$.  Under the assumption that $K$ is pointed, 
$S$ is non-empty by Problem 7.15 of \cite{Sturmfels2004}.
Let $\bar S = Q\setminus S = \mbox{{\em non-saturation
points} of } Q$.

Now, consider minimal points of $S$ with respect to $S$ or $Q$.
We call $\Ba\in S$ an {\em $S$-minimal} (or a {\em $Q$-minimal}\/)
% a {\em $\Qsat$-minimal}, resp.) 
saturation point 
if there exists no other 
$\Bb\in S$, $\Bb\neq \Ba$, such that $\Ba- \Bb \in S$  
(or $Q$).
%($Q$, $\Qsat$, resp.).  
Let $\min(S;S)$ be the set of $S$-minimal saturation points
and 
$\min(S;Q)$ be the set of $Q$-minimal saturation points.
% , and 
% $\min(S;Q_{\rm sat})$ be the set of $Q_{\rm sat}$-minimal saturation points.

% For each fundamental hole $\By \in H_0$ and for each $\Ba \in A$,
% we define 
% $\bar\lambda(\By, \Ba)= \min \{ \lambda \in \Z \mid \By +
% \lambda \Ba \in Q\} $, if the set is non-empty. 
% Otherwise we define $\bar\lambda(\By, \Ba)=\infty$.
% Similarly for each element $\Bb\in B$ in the Hilbert basis of $K$ and 
% for each $\Ba \in A$, we define
% $\bar\mu(\Bb, \Ba)= \min \{ \lambda \in \Z \mid \Bb +
% \lambda \Ba \in Q\} $, if the set is non-empty. 
% Otherwise we define $\bar\lambda(\Bb, \Ba)=\infty$.

The following is a list of some notation:
\begin{eqnarray*}
Q&=&\{\lambda_1 \Ba_1 + \cdots + \lambda_n \Ba_n: 
\lambda_1, \cdots, \lambda_n \in \N = \{0,1,\dots\}\} \\
K&=& \{\lambda_1 \Ba_1 + \cdots + \lambda_n \Ba_n: 
\lambda_1, \cdots, \lambda_n \in \R_+ \}\\
%L&=&A {\mathbb Z}^n =\{\lambda_1 \Ba_1 + \cdots + \lambda_n \Ba_n: 
%\lambda_1, \cdots, \lambda_n \in \Z \} \\
\Qsat&=&K \cap \Z^d = \mbox{saturation of $A$} \supset Q \\
H &=& \Qsat \setminus Q = \mbox{holes in } \Qsat\\
H_0 &=& \{ \Ba \in \Qsat : \Qsat \cap (\Ba + (-Q)) = \{ \Ba \}, \, \Ba\neq 0
\}\\
&=&  \mbox{fundamental holes in } \Qsat\\
S &=& \{ \Ba \in Q : \Ba + \Qsat \subset Q \}= \mbox{saturation
  points of } Q\\
\bar S &=& Q\setminus S = \mbox{non-saturation
  points of } Q\\
\min(S;S) &=& \{\Ba \in S: (\Ba + (-(S\cup\{0\})))\cap S=\{ \Ba\}\}\\
   &=& \mbox{minimal $S$-saturation  points of } Q\\
\min(S;Q) &=& \{\Ba \in S: (\Ba + (-Q))\cap S=\{ \Ba\}\}\\
   &=& \mbox{minimal $Q$-saturation  points of } Q \\
% \min(S;\Qsat) &=& \{\Ba \in S: (\Ba + (-\Qsat))\cap S=\{ \Ba\}\}\\
%    &=& \mbox{minimal $\Qsat$-saturation  points of } Q\\
% \bar\lambda(\By, \Ba)&=& \min \{ \lambda \in \Z \mid \By +\lambda \Ba
% \in Q\}, \quad \By\in H_0, \ \Ba\in A \\
% \bar\mu(\Bb, \Ba)&=& \min \{ \lambda \in \Z \mid \Bb +
% \lambda \Ba \in Q\} , \quad \Bb\in B, \ \Ba\in A
\end{eqnarray*}
%\vskip 0.2cm

Using these definitions and notation, \cite{takemura-yoshida2006} showed 
several necessary and sufficient conditions for $H$'s finiteness,
which are summarized in the  following theorem:

\begin{thm}[\cite{takemura-yoshida2006}]
\label{thm:previous-paper}
% Let $\Ba_1,\ldots,\Ba_n \in \R^d$.
% We assume the lattice $L$ generated by $\Ba_1,\ldots,\Ba_n$ equals to $\Z^d$ 
% without loss of generality.
% Let $H_0 =\{ \By_1, \dots, \By_m\}$ be the set of fundamental holes of 
% the semigroup $Q$ generated by $\Ba_1,\ldots,\Ba_n$.  
Under the assumptions and the definitions above 
the following statements are equivalent:
\begin{enumerate}
\setlength{\itemsep}{0pt}
\item\label{equiv1} $\min(S;S)$ is finite. 
\item\label{equiv2} $\cone(S)$ is a closed rational polyhedral cone.
\item\label{equiv3} There is some $s\in S$ on every extreme ray of $K$.
\item\label{equiv4} $H$ is finite.
% \item $\bar\lambda(\By, \Ba) < \infty$, \quad $\forall \By\in H_0$, $\forall
%   \Ba\in A$.
% \item $\bar\mu(\Bb, \Ba) < \infty$, \quad $\forall \Bb \in B$,  $\forall \Ba
%   \in A$.
\item\label{equiv5} $\bar S$ is finite.
\end{enumerate}
\end{thm}
%[some detailed summary of the previous paper here or maybe in
%appendix, fundamental holes and other concepts have to be explained]

In this paper, we further investigate saturation points and holes in a 
semigroup with the given polyhedral cone $K$ and we study how 
saturation points are distributed in each face of the fixed polyhedral cone 
$K$.  In Section \ref{prelim}, we will define {\em almost saturated} faces and
{\em nowhere saturated} faces of $K$.  Then, we will extend the results
in \cite{takemura-yoshida2006} in terms of almost saturated faces and
nowhere saturated faces of $K$, and we will show some 
preliminary results on almost saturated faces and
nowhere saturated faces of $K$.  In Section \ref{almostSat},
we will give the necessary and sufficient conditions for a face of $K$ to be
almost saturated or to be nowhere saturated. 
The results in Section \ref{construction} show that one can construct a 
semigroup $Q$  
from any antichain of faces of any given cone $K$ so that 
the faces in the antichain are {\em minimal almost saturated} in $Q$.
In Section \ref{examples} we will apply our theorems to a small example
and to a more complicated example of 
$2 \times 2 \times 2 \times 2$ tables with three 
$2$-marginals and a $3$-marginal as the simplicial complex on $4$ nodes 
$[12][13][14][234]$ with levels of $2$ on each node.

\section{Definitions and preliminary results}\label{prelim}

% Let $K$ be a pointed rational polyhedral cone with non-empty interior
% in $\rr^d$.  $K$ is arbitrary but fixed throughout this paper.  
% We denote the set of integer points of $K$ by
% \[
% Q_K = K \cap \zz^d.
% \]

Let $\cal F$ denote
the face poset
%the face lattice (Section 8.6 of \cite{schrijver1986})
of $K$.  For each proper face $F \in {\cal F}$, $F \subsetneq K$,
there exists a supporting hyperplane
\[
L_F = \{ \Bu \in \rr^d \mid \Bc_F \cdot \Bu=0\}
\]
such that $\Bc_F \cdot \Bu \ge 0$, $\forall \Bu\in K$, and
\[
F = L_F \cap K .
\]
For each proper face $F$, we choose and fix  $\Bc_F$ throughout this
paper.  %For notational consistency we define $\Bc_K = 0$.  
$\relint(F)$ denotes the relative interior of $F$.

Let $A_F = A \cap F$ denote the set of elements of $A$ in the face $F$.
% we consider various $A$'s with
% $K=\cone(A)$.
Similarly define $Q_F = Q(A) \cap F$.  Then 
\[
Q_F = \sg(A_F)
\]
is the commutative semigroup generated by $\Ba \in A_F$.
As a particular convenient element in $Q(A)\cap \relint(F)$, we often make use
of the following element
\begin{equation}
\label{eq:a-star-f}
\Ba_F^* = \sum_{\Ba_i \in A_F} \Ba_i .
\end{equation}

Elements of
the Hilbert basis $B$ in $F$ is denoted by  
\[ B_F = B \cap F .
\]
$B_F$ is the unique minimal Hilbert basis  for $F$.
It can be easily shown that 
every element in the Hilbert basis $B$ not belonging to $Q$ is a fundamental hole:
\begin{equation}
\label{eq:hilbert-not-in-Q}
B \setminus Q \subset H_0 .
\end{equation}
%$B_{\rm ray} \subset B$ denotes the set of elements of $B$ which are
%extreme rays of $K$.
%Note that $\Bb \in B_{\rm ray}$
%is the integral point closest to the origin on the ray of $K$ spanned
%by $\Bb$.

Now we give a key definition for this paper.

\begin{defn}
\label{def:face-saturation}
  We call a face $F$ {\em almost saturated} if there exists a
  saturation point of $Q=Q(A)$ on F. Otherwise (i.e.\ if no point of
  $F$ is a saturation point) we call $F$ {\em
    nowhere saturated}.
\end{defn}

For the one-dimensional faces (i.e.\ extreme rays) of $K$ this
definition corresponds to Condition \ref{equiv3} of Theorem 
\ref{thm:previous-paper}.
The basic fact on the existence of saturation point
(Problem 7.15 of \cite{Sturmfels2004}) is that $K$ itself
is always almost saturated.  However generally we do not know whether any
other face of $K$ is almost saturated or not.  
% In Corollary 10 of 
% \cite{takemura-yoshida2006} we established the following proposition.
% \begin{prop}
% The set of holes $H$ is finite if and only if all 1-dimensional
% faces of $K$ are almost saturated.
% \end{prop}
Therefore the important question is to ask which faces of $K$ are
%almost saturated in the case that $H$  is infinite.
almost saturated if $H$  is infinite.

Now from the definition the following lemma is obvious.

\begin{lemma}
\label{lem:inclusion}
Let $F,G$ be two faces of $K$ with $G\subset F$.
If $G$ is almost saturated, then $F$ is almost saturated.
Alternatively, if $F$ is nowhere saturated, then $G$ is nowhere
saturated.
\end{lemma}

This lemma shows that if there are infinitely many
holes, we can describe how the holes are distributed in terms of the
set of minimal almost saturated faces or the set of maximal nowhere
saturated faces.  Here ``minimal'' and ``maximal'' refer to the
partial order of the face poset in terms of inclusion of faces.

In the above lemma suppose that $\Bb\in G$ is a saturation point.
Then for any $\Ba \in \relint(F)$, $\Bb + \Ba \in \relint(F)$ and $\Bb
+ \Ba$ is a saturation point.  Therefore if $F$ is almost saturated,
then there always exists a saturation point in $\relint(F)$.
Alternatively, $F$ is nowhere saturated if no point of $\relint(F)$ is
a saturation point.  We can summarize this fact as follows.
\[
F \cap S = \emptyset \quad \Leftrightarrow \quad \relint(F)\cap S =
\emptyset.
\]

We end this section with some 2-dimensional examples to illustrate our
definitions.  The commutative semigroups of the examples share a
common saturation but the distribution of the holes are different.  In
our examples we will write $A$ as a $d\times n$ integral matrix, so that
$\Ba_1, \dots, \Ba_n$ are columns vectors of $A$.

\begin{figure}[ht]
\begin{center}
     \includegraphics[scale=.4]{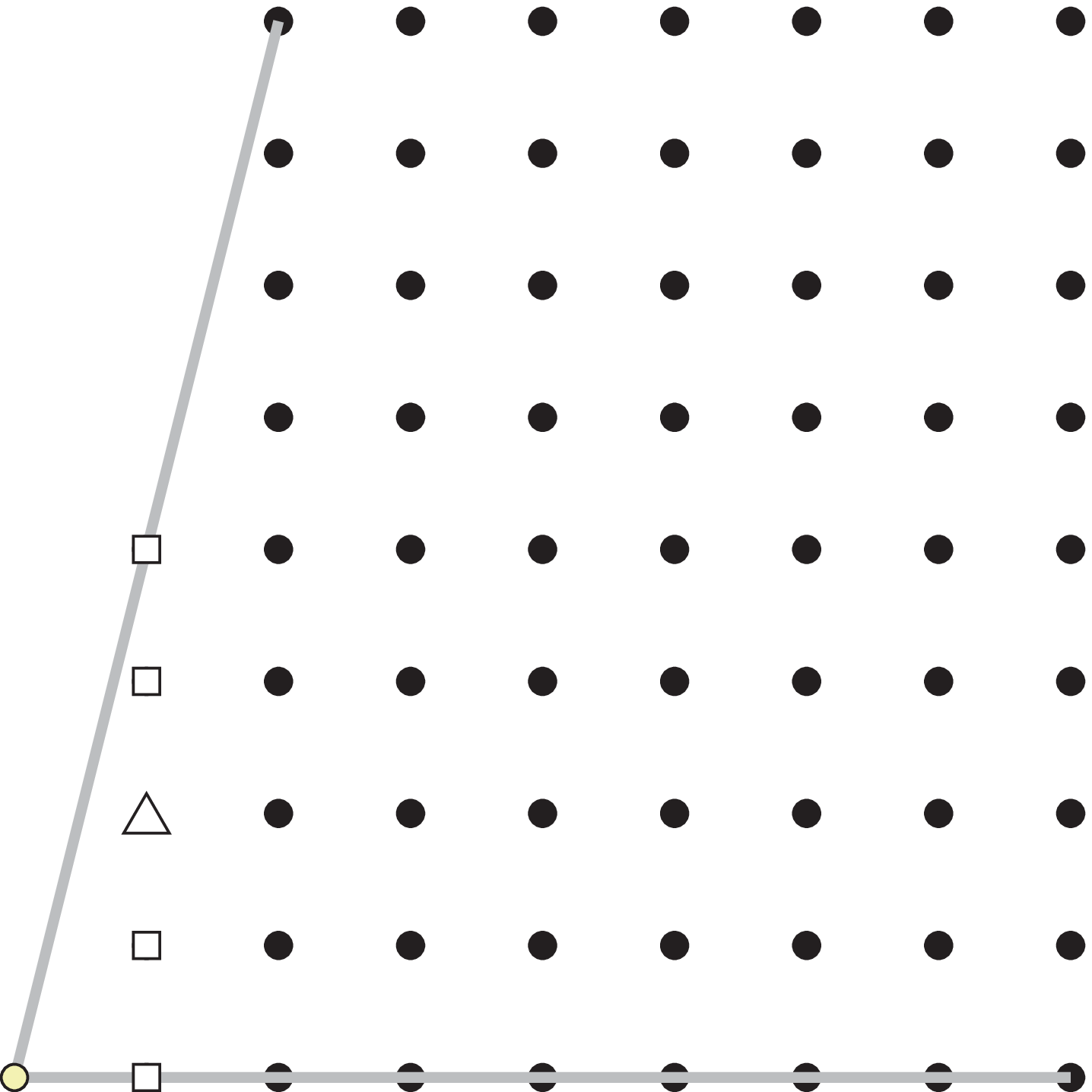}
 \end{center}
\caption{White circles represent nonsaturation points, a triangle represents a hole, 
white squares represent $S$-minimal saturation points, and black circles 
represent non $S$-minimal saturation points  
in the semigroup in Example \ref{example1}.}
\end{figure}

\begin{figure}[ht]
\begin{center}
     \includegraphics[scale=.4]{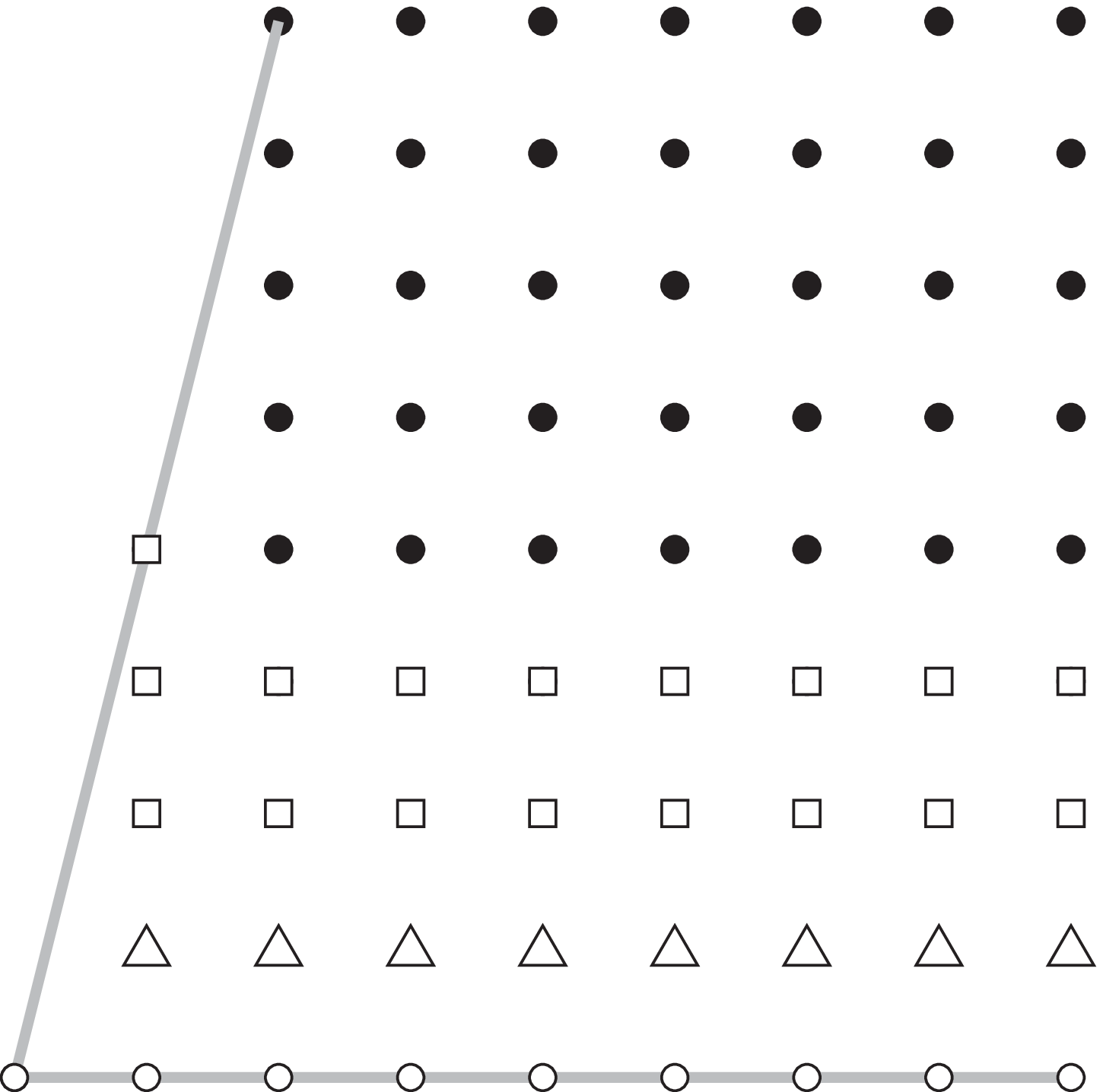}
 \end{center}
\caption{White circles represent nonsaturation points, triangles represent holes, 
white squares represent $S$-minimal saturation points, and black circles 
represent non $S$-minimal saturation points  
in the semigroup in Example \ref{example2}.}
\end{figure}

\begin{figure}[ht]
\begin{center}
     \includegraphics[scale=.4]{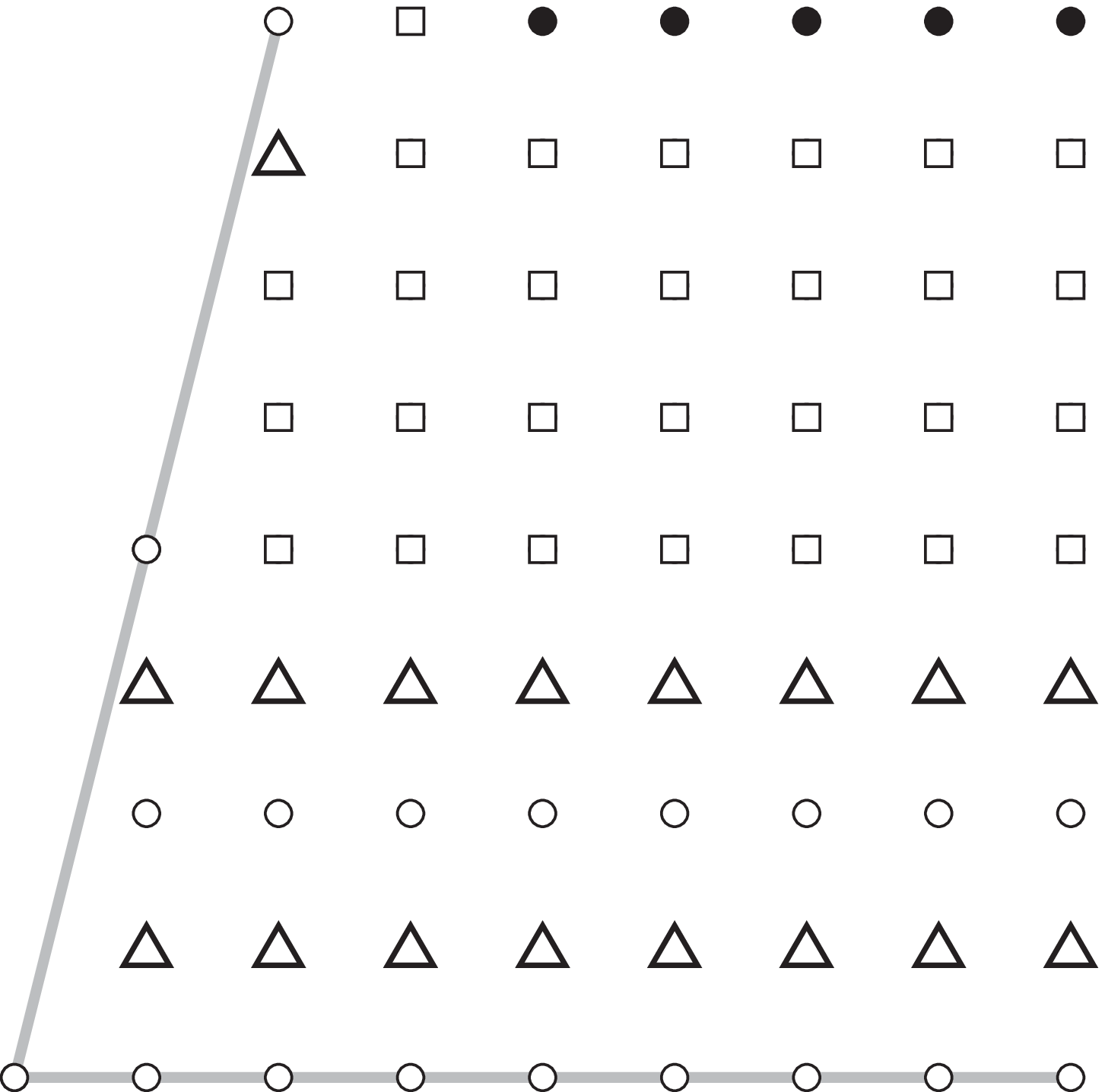}
 \end{center}
\caption{White circles represent nonsaturation points, triangles represent holes, 
white squares represent $S$-minimal saturation points, and black circles 
represent non $S$-minimal saturation points  
in the semigroup in Example \ref{example3}.}
\end{figure}

\begin{ex}\label{example1}
Let $A$ be an integral matrix such that
\[
A = \left(\begin{array}{cccc}
1 & 1 & 1 & 1\\
0 & 1 & 3 & 4\\
\end{array}\right).
\]
The cone $K$ is defined by $K = \cone(A) = \{(x_1, x_2)^t: -4x_1 + x_2 \leq 0, 
x_1, x_2 \in \R_+\}$.
The set of holes $H$ consists of only one element $\{(1, 2)^t\}$.
$\bar S = \{(0, 0)^t\}$. $\min(S;S) = \{(1,0)^t,\, (1,1)^t,\,
(1, 3)^t, \, (1, 4)^t\}$
Thus, $H$, $\bar S$, and $\min(S;S)$ are all finite.
Note that one-dimensional faces $F_1 = \{(x_1, 0)^t: x_1 \in \R_+\}$ and
$F_2 = \{(x_1, x_2)^t: -4x_1 + x_2 = 0, x_1, x_2 \in \R_+\}$ contain
saturation points in each relative interior
(i.e. a saturation point $(2, 0)^t \in 
\relint(F_1)$ and a saturation point $(1, 4)^t \in \relint(F_2)$).
Since $F_1 \subset K$, by Lemma \ref{lem:inclusion}, $K$ is almost
saturated and since $K$ contains saturation points $(2, 0)^t$ and
$(1, 4)^t$, $K$ is indeed almost saturated.
\end{ex}
\begin{ex}\label{example2}
Let $A$ be an integral matrix such that
\[
A = \left(\begin{array}{cccc}
1 & 1 & 1 & 1\\
0 & 2 & 3 & 4\\
\end{array}\right).
\]
The cone $K $ is defined by $K = \cone(A) = \{(x_1, x_2)^t: -4x_1 + x_2 \leq 0, 
x_1, x_2 \in \R_+\}$ which is the same as Example \ref{example1}.
The set of holes $H$ consists of elements $\{(k, 1): k \in \Z , \, 
k \geq 1 \}$.
$\bar S = \{(i, 0)^t : i  \in \Z , \, i\geq 0\}$.
and 
$\min(S;S)=\{(k, j)^t: k \in \Z , \, k \geq 1, \, 
2 \leq j \leq 3 \} \cup \{(1, 4)\}$.
Thus, $H$, $\bar S$, and $\min(S;S)$ are all infinite.
However, $\min(S;Q) = \{(1, 2)^t, \, (1, 3)^t, \, (1, 4)^t\}$ is finite.
Note that a one-dimensional face
$F_2 = \{(x_1, x_2)^t: -4x_1 + x_2 = 0, x_1, x_2 \in \R_+\}$ contains
a saturation point $(1, 4)^t \in \relint(F_2)$ in its relative interior,
so $F_2$ is almost saturated.  However,
another one-dimensional face $F_1 = \{(x_1, 0)^t: x_1 \in \R_+\}$
does not contains any saturation points in its relative interior, so $F_1$
is nowhere saturated.  Since $\bar S = Q_{F_1} = Q \cap F_1$,
$Q_{F_1}$ does not contain
any saturation points and also note that since $F_2$ is almost saturated,
$K$ is almost saturated ($K$ contains a saturation point $(1, 4)^t$ in 
its relative interior).
\end{ex}
\begin{ex}\label{example3}
Let $A$ be an integral matrix such that
\[
A = \left(\begin{array}{cccc}
1 & 1 & 2 & 1\\
0 & 2 & 5 & 4\\
\end{array}\right).
\]
The cone $K$ is defined by $K = \cone(A) = \{(x_1, x_2)^t: -4x_1 + x_2 \leq 0, 
x_1, x_2 \in \R_+\}$ which is the same as Example \ref{example1} and Example 
\ref{example2}.
The set of holes $H$ consists of elements $\{(k, 1)^t: k \in \Z , \, 
k \geq 1 \} \cup \{(k, 4k-1)^t: k \in \Z , \, k \geq 1 \}\cup 
 \{(k, 3)^t: k \in \Z , \, k \geq 1 \}$.
$\bar S = \{(i, 0)^t : i  \in \Z , \, i\geq 0\} \cup \{(i, 2)^t : i  \in \Z , 
\, i\geq 1\} \cup 
\{(i, 4i)^t : i  \in \Z , \, i\geq 0\}$
and 
$\min(S;S)=\{(k, j)^t: k \in \Z , \, k \geq 2, \, 
4 \leq j \leq 6 \} \cup \{(3, 8)^t\}\cup
\{(k, j)^t: k \in \Z , \, k \geq 3, \, 
 j = 7 \} \cup \{(k, 4k-i)^t: k, \, i  \in \Z , \, k \geq 3, \, 
i = 2, \, 3 \}$.  Thus, $H$, $\bar S$, and $\min(S;S)$ are all infinite.
However, $\min(S;Q) = \{(2, 4)^t, \, (2, 5)^t, \, (2, 6)^t\}$ is finite.
Note that one dimensional faces  $F_1 = \{(x_1, 0)^t: x_1 \in \R_+\}$ and
$F_2 = \{(x_1, x_2)^t: -4x_1 + x_2 = 0, x_1, x_2 \in \R_+\}$
do not contain any saturation points in its relative interior, so $F_1$
and $F_2$
are nowhere saturated.  
However, since $K$ contains a saturation point $(2, 4)^t$, so $K$ is almost 
saturated.
\end{ex}

%[put two dimensional example here]

\section{Condition for almost saturation of a face}\label{almostSat}

In this section we consider the problem of determining which face of
$K=\cone(A)$ is almost saturated for a given $A$.  We use the
following fact proved in Lemma 4.1 of \cite{takemura-yoshida2006}.

\begin{lemma}
\label{lem:3.1}
Suppose that the semigroup $Q=Q(A)$ is not saturated. 
$\Ba\in Q$ is a saturation point if
and only if $\Ba + \By \in Q$ for all fundamental holes $\By$.
\end{lemma}

% Let $F$ be a face of $K$. We denote the elements of the semigroup $Q$
% on the face $F$ by
% \[
% Q_F = Q \cap F.
% \]
% $Q_F$ is generated by $\Ba_i$'s in $F$.
For a face $F$, we denote the parallel shift of $Q_F$ by $\By \in
K$ by
\[
\By + Q_F = \{ \By + \Ba \mid \Ba\in Q\cap F\}.
\]
We will consider the case that $\By$ is a hole.
We now prove the following fact. 
% [The proof is a bit difficult.  We
% need to prepare some facts more carefully.]
\begin{prop}
\label{prop:basic}
  A face $F$ is nowhere saturated if and only if there exists some
  fundamental hole $\By \in H_0$ such that all points of $\By + Q_F$
  are holes, i.e.\ $\By+Q_F \subset H$.
\end{prop}

\begin{proof}
Suppose that $\By + Q_F \subset H$ for some fundamental hole $\By\in
H_0$.  Then for $\By + \Ba \not\in Q$ for every $\Ba \in Q_F$.
By Lemma \ref{lem:3.1} $\Ba$ is not a
saturation point.  Therefore $F$ is nowhere saturated.
% by Lemma \ref{lem:trivial}.

Conversely suppose that $F$ is nowhere saturated.  
%By Lemma %\ref{lem:trivial} 
Then no $\Ba\in Q_F$ is a saturation point. 
% Also by Lemma
% \ref{lem:inclusion} there is no saturation point on the boundary of
% $F$ either. Therefore no point of $Q_F$ is a saturation point.
% (Maybe we should state this earlier).  
By Lemma \ref{lem:3.1}
for each $\Ba \in Q_F$ there exists a fundamental hole  $\By \in H_0$
such that $\Ba + \By \in H$. 
% Since there exist only a finite number of
% fundamental holes, there exists a fundamental hole $\By$ such that
% $\By + \Ba \in H$ for infinitely many $\Ba \in \relint(Q)\cap F$.
Consider $\Ba_F^*$ in (\ref{eq:a-star-f}).  $\Ba_F^* \in \relint(F)\cap Q$.
Consider $\Ba_F^*, 2\Ba_F^*, 3\Ba_F^*,
\dots$.  These are not saturation points of $Q$.  Therefore 
for each $k \Ba_F^*$, $k\in \nn$, there exists a fundamental hole
$\By=\By_k$ such that $k \Ba_F^* + \By \in H$. Since there are only finite
number of fundamental holes, there exists a fundamental hole $\By_0\in
H_0$ such that $k\Ba_F^* + \By_0 \in H$ for infinitely many $k$.  

Now each
$\Ba\in Q_F$ is a non-negative integer combination of
$\Ba_i$'s in $Q_F$:
\[
\Ba = \sum_{\Ba_i \in A_F} x_i \Ba_i, \quad x_i \in \nn.
\]
Choose  $m \ge \max\{ x_i \mid \Ba_i \in A_F\}$ such that $m \Ba_F^* +
\By_0 \in H$.
% Therefore for each $\Ba \in \relint(F)\cap Q$
% there exists a large enough $k$ such that
Then
\[
m \Ba_F^* - \Ba = \sum_{\Ba_i \in A_F} (m - x_i) \Ba_i
\in Q_F.
\]
% In this case we have $\Ba + \By_0 \in H$.  Therefore $\By_0 + Q_F
% \subset H$.
Now if $\Ba + \By_0 \in Q$, then
\[
m\Ba_F^* + \By_0 = \Ba + \By_0 + (m \Ba_F^* - \Ba) \in Q, 
\]
which is a contradiction.  Therefore $\Ba + \By_0 \in H$ for every
$\Ba\in Q_F$.
\end{proof}

We now state the following main result of this section. 
\begin{thm}
\label{thm:feasible}
A face $F$ is almost saturated if and only if every fundamental hole
$\By\in H_0$ can be written as 
\begin{equation}
\label{eq:feasible}
\By = x_1 \Ba_1 + \dots + x_n \Ba_n,  \quad x_j \in \zz,  \forall j,\
\ \text{and} \ \ 
%q_j \le 0 \ \text{for}\  \Ba_j \in F, \ \ 
x_j \ge 0 \ \text{for}\  \Ba_j \not\in F.
\end{equation}
Equivalently, $F$ is nowhere saturated if and only if
(\ref{eq:feasible}) does not have a feasible solution for some
fundamental hole $\By$.
\end{thm}

\begin{proof}
Suppose that $F$ is almost saturated. Then by Proposition 
\ref{prop:basic},
$\By + Q_F \not\subset H$ for every fundamental hole $\By$.
Therefore for every $\By$
there exists $\Ba \in Q_F$ such that
$\By + \Ba \in Q$, i.e. $\Bb = \By +\Ba \in Q$.  But then
$\By = \Bb - \Ba$  is a feasible solution to (\ref{eq:feasible}).
Therefore, every fundamental hole $\By$ can be written in the form in 
\eqref{eq:feasible}.

Conversely suppose that 
(\ref{eq:feasible}) has a feasible solution for every fundamental hole
$\By$.  Then
\[
\By + \sum_{\Ba_i \in F} |x_i| \Ba_i = 
\sum_{\Ba_i \in F} (x_i + |x_i|) \Ba_i + 
\sum_{\Ba_i \not\in F} x_i
\Ba_i \in Q.
\]
Since $\sum_{\Ba_i \in F} |x_i| \Ba_i \in Q_F$ and since
$\By + \sum_{\Ba_i \in F} |x_i| \Ba_i \in Q$, 
$F$ is almost saturated by Proposition \ref{prop:basic}.
\end{proof}

Theorem \ref{thm:feasible} is stated in terms of the fundamental
holes.  For applications it is more convenient to have a condition for
elements in the Hilbert basis $B$ of $K$. % the cone $\cone(A)$.  
We have the following theorem. 
% [actually
% this is a conjecture at this point, but I am pretty sure that it holds.]

\begin{thm}
\label{thm:feasible-hilbert}
A face $F$ is nowhere saturated if and only if for some element $\Bb$ of the
Hilbert basis $B$
\begin{equation}
\label{eq:feasible-hilbert}
\Bb = x_1 \Ba_1 + \dots + x_n \Ba_n,  \quad x_j \in \zz, \ \forall j,\ 
\ \  \text{and}\ \ 
%q_j \le 0 \ \text{for}\  \Ba_j \in F, \ 
x_j \ge 0 \ \text{for}\  \Ba_j \not\in F.
\end{equation}
does not have a feasible solution. 
\end{thm}

\begin{proof}
Suppose that (\ref{eq:feasible-hilbert}) does not have a feasible
solution for some $\Bb \in B$.  Obviously this $\Bb$ can not belong
to $Q$.  Then by (\ref{eq:hilbert-not-in-Q}) $\Bb$ is a fundamental
hole and $F$ is nowhere saturated  by Theorem \ref{thm:feasible}.

Conversely suppose that 
(\ref{eq:feasible-hilbert}) has a feasible solution for every element
$\Bb \in B$.  Now every fundamental hole $\By\in H_0$ can be
written as a non-negative integral combination of $\Bb$'s.  Then forming
the same non-negative integral combination of the feasible solutions of
$\Bb$'s, we obtain a feasible solution to 
(\ref{eq:feasible}) for every fundamental hole $\By \in
H_0$. Therefore by Theorem \ref{thm:feasible} $F$ is almost saturated.
\end{proof}

Finally we prove the following theorem.
\begin{thm}
\label{thm:minimal-q}
A face $F$ is nowhere saturated if and only if there exists a $Q$-minimal
  saturation point $\By$  such that every point of  $\By+Q_F$
  is an $S$-minimal  saturation point.
\end{thm}

\begin{proof}
The set of $Q$-minimal saturation
points $\min(S;Q)$ is always finite by Proposition 4.4 of
\cite{takemura-yoshida2006}.
Suppose that $F$ is nowhere saturated.  Then $\min(S;Q)\cap F =\emptyset$.
In Section 1, we chose the normal vector $\Bc_F$  for $F$. Choose a
point $\By_0$ from $\min(S;Q)$ which minimizes $\Bc_F \cdot \By$:
\[
\Bc_F \cdot \By_0 = \min_{\By\in \min(S;Q)} \Bc_F \cdot \By > 0.
\]
Because $\By_0$ is a saturation point, all points on $\By_0 + Q_F$
are saturation points.  Suppose that some $\Ba\in \By_0 + Q_F$ is not
an $S$-minimal saturation point.  Then $\Ba$ can be written as a sum
of non-zero saturation points 
\[
\Ba = \Ba_1 + \Ba_2, \qquad \Ba_1, \Ba_2\in S
\]
and $\Bc_F \cdot \Ba = \Bc_F \cdot \Ba_1 + \Bc_F \cdot \Ba_2$. Since
$\Ba_1, \Ba_2$ do not belong to $Q_F$,  both 
$\Bc_F \cdot \Ba_1$ and $\Bc_F \cdot \Ba_2$ are positive.  In
particular  $\Bc_F \cdot \Ba > \Bc_F \cdot \Ba_1 > 0$.  However
$\Bc_F \cdot \Ba = \Bc_F \cdot \By_0$ and
\[
\Bc_F \cdot \By_0 > \Bc_F \cdot \Ba_1 > 0.
\]
Now $\Ba_1$ can be written as $\Ba_1 = \By_1 + \tilde \By$, $\By_1 \in
\min(S;Q)$, $\tilde \By\in Q$, and therefore  $\Bc_F \cdot \Ba_1 \ge
\Bc_F\cdot \By_1$.  Then $\Bc_F \cdot \By_0 > \Bc_F \cdot \By_1 > 0$,
but this contradicts our choice of $\By_0$.  Therefore 
$\By_0 + Q_F\subset \min(S;S)$.

Conversely suppose that $F$ is almost saturated. Then there exists a
saturation point $\Ba$ on $Q_F$.  Let $\By$ be any $Q$-minimal
saturation point and consider $\By + \Ba$, which is a sum of two
non-zero saturation points and therefore not $S$-minimal.  But 
$\By + \Ba \in \By + Q_F$.  Therefore $\By + Q_F \not\subset
\min(S;S)$. 
\end{proof}

\section{Construction of a semigroup with arbitrary configuration of
  almost saturated faces}\label{construction}

The set of minimal almost saturated faces form an antichain (Section
3.1 of \cite{stanley-vol1}), i.e.\ there is no inclusion relation
among minimal almost saturated faces.  Then a natural question to
ask is whether for any given antichain of faces we can construct
$Q=Q(A)$ such that the set of minimal almost saturated faces of $Q$
coincides with the given antichain.  In this section we show that it
is always possible by explicitly constructing $A$ for a given
antichain.

Let $\{ F_1, \dots, F_M\}$ be an antichain
of faces.  
Write
\[
\bar B = B \cup \{ 0 \}, 
\]
which is the Hilbert basis with the origin added.
For $F_i$, $i=1,\dots,M$, let
\[
\Bb_{F_i}^* = \sum_{\Bb \in B_{F_i}} \Bb  \ \in \relint(F_i).
\]

\comment{
\begin{thm}
\label{thm:configuration}
For any given antichain of faces $\{F_1, \dots, F_M\}$ let
\[
A = \{ C \Bb \mid \Bb \in B_{\rm ray}\} \cup \left(
\{ C \Bb_{F_1}^*, \dots, C \Bb_{F_M}^* \} +  \bar B\right)
\]
Then for sufficiently large $C$ we have
{\rm i)} $\cone(A)=K$, {\rm ii)} the saturation of $Q(A)$ is $Q_K$, 
{\rm iii)} the minimal almost saturated faces for $Q(A)$ are
$F_1, \dots, F_M$.
\end{thm}
\begin{proof}
$\cone(A)=K$ because $A$ contains 
$\{ C \Bb \mid \Bb \in B_{\rm ray}\}$, which is the set of multiples
of all the rays of $K$. The saturation of $Q(A)$ is $Q_K$ because 
$A$ contains 
$C \Bb_{F_1}^*$ and $C \Bb_{F_1}^* + \Bb$ for all $\Bb\in B$ and
each  $\Bb\in B$ can be written as
\[
 (C \Bb_{F_1}^* + \Bb) -  C \Bb_{F_1}^* .
\]
Now we need to show the 
we can choose large enough $C$ such that 
$G$ is nowhere saturated  for all $G\in {\cal G}$, where
\[
{\cal G} =\{ G \mid G \not\supset F_i,  \quad \forall i=1,\dots,M\}.
\]
For each $G\in {\cal G}$ consider $\Bc_G$
chosen in Section 1. $G$ is a proper face and there exists $\Bb\in B$
such that $\Bb \not \in G$. For each $G$ choose $\Bb_G \in  B \setminus G$
arbitrarily.   Let 
\begin{equation}
\label{eq:large-C}
C > \frac{\max_{G \in {\cal G}} \Bc_G \cdot \Bb_G}
 {\min_{G\in \cal G} \min_{\Bb\in B\setminus G} \Bc_G \cdot\Bb}.
\end{equation}
Each element $\Ba$ of $Q(A)$ not on $G$ has at least one $\Bb\in B\setminus
G$ with positive weight.  Then for the above choice of $C$ we have
$\Bc_G\cdot \Ba > \min_{\Bb\in B\setminus G} \Bc_G'\Bb$.
Therefore $\Ba$ does not belong to  the set
\[
\Bb_G + Q_G.
\]
Therefore points on $\Bb_G + Q_G$ are all holes. Since this holds
simultaneously for all $G\in {\cal G}$, we see that with the above
choice of $C$
every $G\in {\cal G}$ is nowhere saturated.
Finally we need to show that each $F_i$ is  almost
saturated.  The argument is similar to the usual proof of the
existence of saturation point for $K$. Having fixed $C$ as in
(\ref{eq:large-C})
we can show that we can choose another large $\bar C$ such that 
$\bar C \times C \Bb_{F_i}^* \in \relint(F_i) $ 
is a saturation point of $Q(A)$.  Therefore  $F_i$ is almost
saturated.
\end{proof}}

\begin{thm}
\label{thm:configuration}
For any given antichain of faces $\{F_1, \dots, F_M\}$ let
\[
A = \{ 2 \Bb \mid \Bb \in B\} \cup \left(
\{ 2 \Bb_{F_1}^*, \dots, 2 \Bb_{F_M}^* \} +  \bar B\right), 
\]
Then we have
{\rm i)} $\cone(A)=K$, {\rm ii)} the saturation of $Q(A)$ is
$K\cap \Z^d$, 
{\rm iii)} the minimal almost saturated faces for $Q(A)$ are
$F_1, \dots, F_M$.
\end{thm}

Note that each element of $A$ is either i) $2\Bb$, $\Bb\in
B$ , ii) $2\Bb_{F_i}^*$,
$i=1,\dots,M$, iii) or of the form $2\Bb_{F_i}^* + \Bb$, $\Bb\in B$.

\begin{proof}  
$\cone(A)=K$ because the Hilbert basis $B$ contains all
the extreme rays of $K$ and $A$ contains  all $2\Bb$, $\Bb\in B$.
The saturation of $Q(A)$ is $K\cap \Z^d$ because $K\cap \Z^d$ is
generated by $B$ and 
each $\Bb\in B$ can be written as
\begin{equation}
\label{eq:bstar-diff}
\Bb = (2 \Bb_{F_1}^* + \Bb) -  2 \Bb_{F_1}^*, 
\end{equation}
where both $2 \Bb_{F_1}^* + \Bb$ and $2 \Bb_{F_1}^*$ belong to $A$.

Now we show that 
every face $G\in {\cal G}$ is nowhere saturated, where
\[
{\cal G} =\{ G \mid G \not\supset F_i,  \quad \forall i=1,\dots,M\}.
\]
For an arbitrary $G\in {\cal G}$ consider $\Bc_G$
chosen in Section 2. $G$ is a proper face and there exists $\Bb\in B$
such that $\Bb \not \in G$. For each $G$ choose $\Bb_G \in  B \setminus G$
such that 
\[
\Bc_G \cdot \Bb_G  = \min_{\Bb\in B\setminus G} \Bc_G \cdot\Bb > 0.
\]
% Let 
% \begin{equation}
% \label{eq:large-C}
% C >  \frac{\max_{G\in {\cal G}}\Bc_G \cdot \Bb_G}
% {\min_{G\in {\cal G}}\Bc_G \cdot \Bb_G} .
% \end{equation}
% {\min_{\Bb\in B\setminus G} \Bc_G \cdot\Bb}.
Now for an arbitrary $G\in {\cal G}$ consider elements  $\Bx$
on $\Bb_G + Q_G$. Then $\Bc_G \cdot \Bx = \Bc_G \cdot \Bb_G$ for all
$\Bx\in \Bb_G + Q_G$.
On the other hand let $\Ba \in Q(A)\setminus G$.  
Since $\Ba$ is a non-negative integral combination of the elements of
$A$, we distinguish two cases: i) $\Ba$ contains some $2\Bb$, $\Bb\not\in G$, 
or ii) $\Ba$ does not contain any $2\Bb$, $\Bb\not\in G$.
In the former case we have
\[
\Bc_G \cdot \Ba \ge 2 \Bc_G \cdot \Bb \ge 2 \Bc_G \cdot \Bb_G
> \Bc_G \cdot \Bb_G
\]
and $\Ba\not\in \Bb_G + Q_G$.  For the latter case $\Ba$ 
has to contain some $2 \Bb_{F_i}^*$. 
%since $\Ba$ is not on $G$.
Furthermore, since $\Ba$ is not on $G$, this $F_i$ is not a subset of $G$ and therefore
$\Bb_{F_i}^*$ contains some $\Bb\not\in G$. Then the same argument as
in the case i) shows that $\Ba\not\in \Bb_G + Q_G$. We have shown that all points of 
$\Bb_G + Q_G$ are holes and hence $G$ is nowhere saturated by
Proposition 
\ref{prop:basic}.

Finally we need to show that each $F_i$ is  almost
saturated.  Let $C=|B|$ denote the number of elements in the Hilbert
basis.  We claim that $2C \Bb_{F_i}^*$ is a saturation point and
therefore $F_i$ is almost saturated.  We are now going to prove this
claim. Every element of $\Qsat=K\cap \Z^d$ is a non-negative integral
combination of elements of $B$.  However $2\Bb$, $\Bb\in B$, are
already contained in $A$.  Therefore 
$\Bx\in Q(A)$ is a saturation point if for every subset $\tilde B$ of
$B$
\[
\Bx + \sum_{\Bb\in \tilde B} b \in Q.
\]
Now by (\ref{eq:bstar-diff})
\[
2C \Bb_{F_i}^* + \sum_{\Bb\in \tilde B} 
\Big((2 \Bb_{F_i}^* + \Bb) -  2 \Bb_{F_i}^* \Big)
= 2(C - |\tilde B|) \Bb_{F_i}^* + \sum_{\Bb\in \tilde B} (2
\Bb_{F_i}^* + \Bb)
\in Q.
\]
Therefore $2C \Bb_{F_i}^*$ is a saturation point.  This proves the theorem.
\end{proof}

%[this proof is still too short.]

\section{Examples}\label{examples}

In this section we will go through some theorems with an example from
Example \ref{example2} and with a defining matrix for $2 \times 2 \times 2 
\times 2$ tables with three 
$2$-marginals and a $3$-marginal as the simplicial complex on $4$ nodes 
$[12][13][14][234]$ with levels of $2$ on each node.

\subsection{$2 \times 4$ matrix from Example \ref{example2}}
Let $A$ be an integral matrix such that
\[
A = \left(\begin{array}{cccc}
1 & 1 & 1 & 1\\
0 & 2 & 3 & 4\\
\end{array}\right).
\]
The cone $K$ is defined by $K= \cone(A) = \{(x_1, x_2)^t: -4x_1 + x_2 \leq 0, 
x_1, x_2 \in \R_+\}$.  One-dimensional faces are $F_1 = \{(x_1, 0)^t: x_1 \in 
\R_+\}$ and $F_2 = \{(x_1, x_2)^t: -4x_1 + x_2 = 0, x_1, x_2 \in \R_+\}$.
The set of holes $H$ consists of elements $\{(k, 1): k \in \Z , \, 
k \geq 1 \}$ and $H_0$ consists of only one element $(1, 1)^t$ since
$\Qsat \cap ((k, 1)^t - Q) = \{(i, 1)^t: i \in \Z, \, 1 \leq i \leq k\}$.
Let $\By = (1, 1)^t$.  Also let $\Ba_1 = (1, 0)^t, \, \Ba_2 = (1, 2)^t, \,
\Ba_2 = (1, 3)^t, \, \Ba_2 = (1, 4)^t$.

Consider $\By + Q_{F_1}$.  Note that $Q_{F_1} = \{(x, 0): x \in \N\}$ so
$\By + Q_{F_1} = \{(x, 1): x \in \N, \, a \geq 1\}$.  Note that 
$\By + Q_{F_1} \subset H$.  Thus, by Proposition 
\ref{prop:basic}, $\By + Q_{F_1} \subset H$ implies $F_1$ is nowhere saturated.

Now we examine whether $F_1$ and $F_2$ are nowhere saturated or almost 
saturated via Theorem \ref{thm:feasible}.
First, we will decide whether $F_1$ is nowhere saturated or almost saturated.
Note that $\Ba_1 \in F_1$, so we set the system of linear equations and 
inequalities such that:
\begin{equation}\label{check1}
\By = x_1 \Ba_1 + x_2\Ba_2 + x_3 \Ba_3 + x_4 \Ba_4, \, x_i \in \Z, \mbox{
for } i = 1, \, \cdots, \, 4,  \, x_j \geq 0, \mbox{
for } j = 2,\, 3,\, 4.
\end{equation}
We count the number of integral solutions in the system \eqref{check1} via 
{\tt LattE} (\cite{lattemanual}) and find out that there is no integral 
solution in the system.  Thus, by Theorem \ref{thm:feasible}, we know that
$F_1$ is nowhere saturated.

Secondly, we will decide whether $F_2$ is nowhere saturated or almost 
saturated.
Note that $\Ba_4 \in F_2$, so we set the system of linear equations and 
inequalities such that:
\begin{equation}\label{check2}
\By = x_1 \Ba_1 + x_2\Ba_2 + x_3 \Ba_3 + x_4 \Ba_4, \, x_i \in \Z, \mbox{
for } i = 1, \, \cdots, \, 4, \, x_j \geq 0, \mbox{
for } j = 1, \, 2,\, 3.
\end{equation}
We count the number of integral solutions in the system \eqref{check2} via 
{\tt LattE} and find out that there are two integral 
solutions in the system, namely $(x_1, x_2, x_3, x_4) = (0, 0, 3, -2)$
and $(x_1, x_2, x_3, x_4) = (0, 1, 1, -1)$.
Thus, by Theorem \ref{thm:feasible}, we know that
$F_2$ is almost saturated.

Now we are going to decide whether $F_1$ and $F_2$ are nowhere saturated or 
almost saturated via Theorem \ref{thm:feasible-hilbert}.
Note that the Hilbert basis $B$ of the cone $K$ consists of $5$ elements 
\[
B = \{\Bb_1 = (1, 0)^t, \Bb_2 = (1, 1)^t, \Bb_3 = (1, 2)^t, \Bb_4 = (1, 3)^t, 
\Bb_5 = (1, 4)^t\}.
\]
Firstly for $F_1$, notice that $\By = \Bb_2$.  Thus since the system 
\eqref{check1} does not have an integral solution so $F_1$ is nowhere 
saturated by Theorem \ref{thm:feasible-hilbert}.

For $F_2$ we set the five systems of equations and inequalities such that:
\begin{eqnarray*}
\Bb_k = x_1 \Ba_1 + x_2\Ba_2 + x_3 \Ba_3 + x_4 \Ba_4, &\\
x_i \in \Z, \mbox{
for } i = 1, \, \cdots, \, 4, \, & x_j \geq 0, \mbox{
for } j = 1, \, 2,\, 3, \, \mbox{
for } k = 1, \, \cdots, 5.
\end{eqnarray*}
Using {\tt LattE} we find out that all systems contain integral solutions,
thus by Theorem \ref{thm:feasible-hilbert}, $F_2$ is almost saturated.

We are going to decide whether  $F_1$ is  nowhere saturated using Theorem
\ref{thm:minimal-q}.  From Example \ref{example2}, we know that
$\min(S;Q) = \{(1, 2)^t, \, (1, 3)^t, \, (1, 4)^t\}$.  We take $(1, 2)^t
\in \min(S;Q)$. Then we notice that $(1, 2)^t + Q_{F_1}$ is contained in
$\min(S;S)$, thus, by Theorem \ref{thm:minimal-q}, $F_1$ is  nowhere saturated.

\subsection{$2\times 2\times 2\times 2$ tables with $2$-marginals and a $3$-marginal}

Now we consider $2 \times 2 \times 2 \times 2$ tables with three 
$2$-marginals and a $3$-marginal as the simplicial complex on $4$ nodes 
$[12][13][14][234]$ with levels of $2$ on each node.
\cite{takemura-yoshida2006} showed that the cardinality of $H$ for the 
semigroup defined by their matrix is infinite.  Thus, we should like to 
investigate which face of the polyhedral cone defined by this matrix
is almost saturated or nowhere saturated by Theorem \ref{thm:feasible-hilbert}.
To compute minimal Hilbert bases of cones, we used 
{\tt normaliz} (\cite{brunskoch}) and to compute each hyperplane 
representation and vertex representation we used 
{\tt CDD} (\cite{fukuda}) and {\tt lrs} (\cite{avis}).
Also we used {\tt 4ti2} (\cite{Hemmecke+Hemmecke+Malkin:2005})
to compute defining matrices.
To count the number of integral solutions in each system, we used
{\tt LattE} (\cite{lattemanual}).

After removing redundant rows (we removed redundant rows using {\tt cddlib})
(\cite{fukuda}), $2 \times 2 \times 2 \times 2$ tables with $2$-marginals and a
$3$-marginal have a $12 \times 16$ defining matrix.  Thus the semigroup is 
generated by $16$ (column) vectors in $\Z^{12}$ such that:
\begin{verbatim}
 1 0 0 0 1 0 0 0 1 0 0 0 1 0 0 0
 0 1 0 0 0 1 0 0 0 1 0 0 0 1 0 0
 0 0 1 0 0 0 1 0 0 0 1 0 0 0 1 0
 0 0 0 1 0 0 0 1 0 0 0 1 0 0 0 1
 1 0 1 0 0 0 0 0 1 0 1 0 0 0 0 0
 0 1 0 1 0 0 0 0 0 1 0 1 0 0 0 0
 1 0 1 0 1 0 1 0 0 0 0 0 0 0 0 0
 0 1 0 1 0 1 0 1 0 0 0 0 0 0 0 0
 1 1 0 0 0 0 0 0 0 0 0 0 0 0 0 0
 0 0 1 1 0 0 0 0 0 0 0 0 0 0 0 0
 0 0 0 0 1 1 0 0 0 0 0 0 0 0 0 0
 0 0 0 0 0 0 0 0 1 1 0 0 0 0 0 0
\end{verbatim}
All of these vectors are extreme rays of the cone (we verified via 
{\tt cddlib}).
The Hilbert basis of the cone generated by these $16$ vectors
consists of these 16 vectors and two additional vectors
\[
\Bb_{17}=(1\ 1\ 1\ 1\ 1\  1\  1\  1\  1\  0\  0\  0)^t, \qquad
\Bb_{18}=(1\  1\  1\  1\  1\  1\  1\  1\  0\  1\  1\  1)^t.
\]
Using {\tt CDD} we computed the system of $48$ linear inequalities for defining
the cone.  Thus, the cone $K$ has $48$ facets.
The results of our experiments are in Table \ref{2x2x2x2Table}.
To enumerate all faces, we used {\tt allfaces\_gmp} from {\tt cddlib}.

\begin{table} 
\begin{center}
\begin{tabular}{|c|c|c|c|} \hline
Dimension & $\#$ of faces & $\#$ of nowhere & $\#$ of almost\\\hline
11 & 48 & 0 & 48\\ \hline
10 & 492 & 0 & 492\\ \hline
9 & 2104 & 0 &  2104 \\ \hline
8 & 4898 & 2 & 4896 \\ \hline
7 & 6956 & 16 & 6940\\ \hline
6 & 6440 & 56 & 6384\\ \hline
5 & 4064 & 112 &3952 \\ \hline
4 & 1796 & 140 & 1656\\ \hline
3 & 560 & 112 & 448 \\ \hline
2 & 120 & 56 & 64 \\ \hline
1 & 16 & 16 & 0 \\ \hline
\end{tabular}
\caption{Faces for $2 \times 2 \times 2 \times 2$ tables with three 
$2$-marginals and a $3$-marginal.  The first column represents the 
dimension of faces, the second column represents the number of faces,
the third column represents the number of nowhere saturated faces, and the 
fourth column represents the number of almost saturated faces.
} 
\label{2x2x2x2Table}
\end{center}
\end{table}

From Table \ref{2x2x2x2Table} we see that $64$ almost saturated $2$
dimensional faces are minimal. We have checked that
all almost saturated faces with dimensions 3 or larger contain at
least one of $64$ almost saturated $2$ dimensional faces.  This
implies that $2$ nowhere saturated $8$ dimensional faces are the
maximal nowhere saturated faces.  Therefore the most important faces
to investigate are $2$ nowhere saturated $8$ dimensional faces and
$64$ almost saturated $2$ dimensional faces.  We give detailed
descriptions of these faces.

% We analyzed $64$ almost saturated $2$ dimensional faces and $2$ nowhere 
% saturated $8$ dimensional faces.

\begin{figure}[ht]
\vskip 0.1in
\begin{center}
\begin{overpic}[width=.50\textwidth]{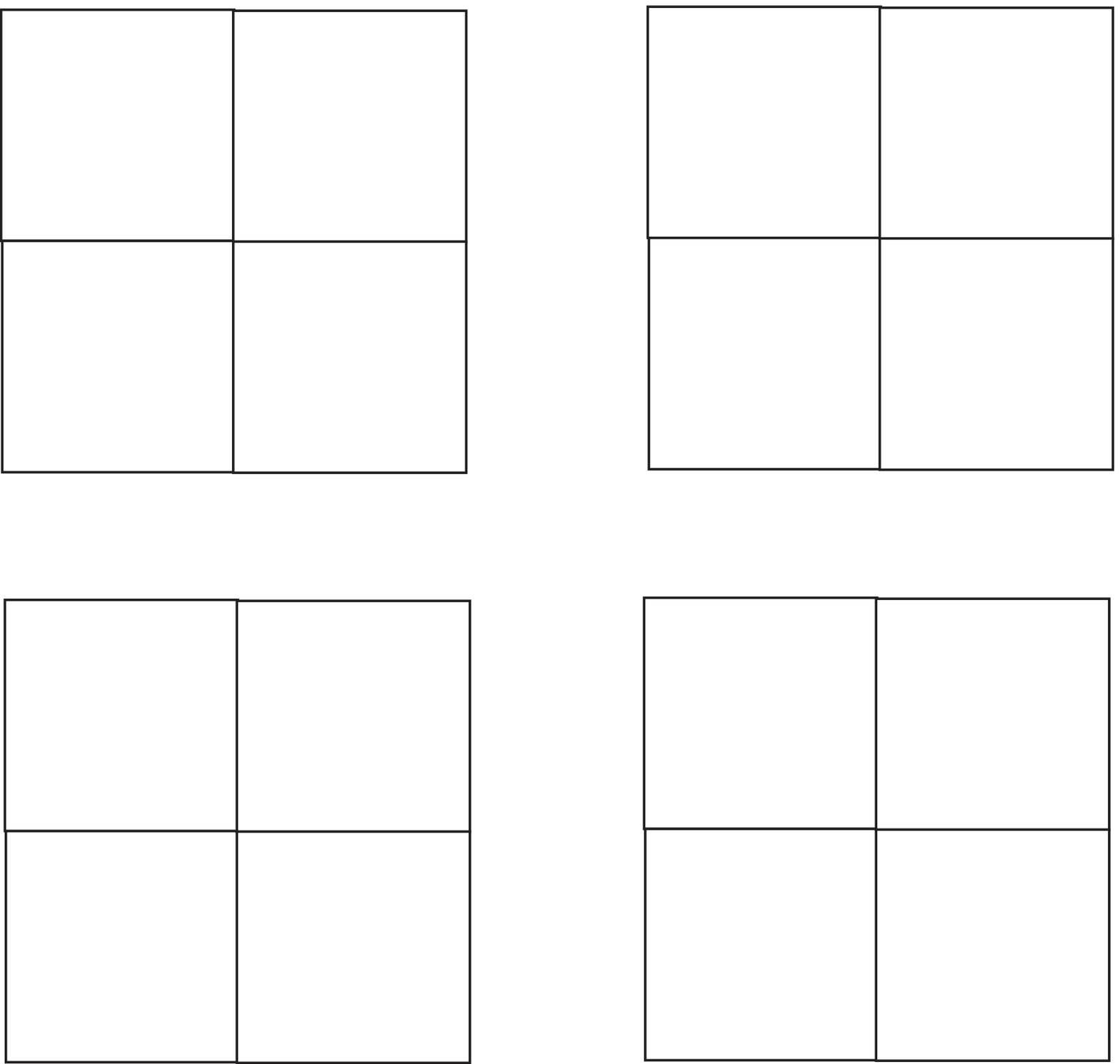}
\put(15,100){$k=1$}
\put(72,100){$k=2$}
\put(5,95){$i=1$}
\put(25,95){$i=2$}
\put(62,95){$i=1$}
\put(82,95){$i=2$}
\put(5,45){$i=1$}
\put(25,45){$i=2$}
\put(62,45){$i=1$}
\put(82,45){$i=2$}

\put(-25,70){$l=1$}
\put(-25,20){$l=2$}
\put(-14,85){$j=1$}
\put(-14,60){$j=2$}
\put(-14,30){$j=1$}
\put(-14,10){$j=2$}
\put(43,85){$j=1$}
\put(43,60){$j=2$}
\put(43,30){$j=1$}
\put(43,10){$j=2$}

\put(10,85){$x_1$}
\put(7,80){$T_{1111}$}
\put(30,85){$x_2$}
\put(27,80){$T_{2111}$}

\put(67,85){$x_3$}
\put(64,80){$T_{1121}$}
\put(87,85){$x_4$}
\put(84,80){$T_{2121}$}

\put(10,65){$x_5$}
\put(7,60){$T_{1211}$}
\put(30,65){$x_6$}
\put(27,60){$T_{2211}$}

\put(67,65){$x_7$}
\put(64,60){$T_{1221}$}
\put(87,65){$x_8$}
\put(84,60){$T_{2221}$}
%%%%%%%%%%%%%%%%%%%%%%%%%%%%%
\put(10,33){$x_9$}
\put(7,28){$T_{1112}$}
\put(30,33){$x_{10}$}
\put(27,28){$T_{2112}$}

\put(67,33){$x_{11}$}
\put(64,28){$T_{1122}$}
\put(87,33){$x_{12}$}
\put(84,28){$T_{2122}$}

\put(10,13){$x_{13}$}
\put(7,8){$T_{1212}$}
\put(30,13){$x_{14}$}
\put(27,8){$T_{2212}$}

\put(67,13){$x_{15}$}
\put(64,8){$T_{1222}$}
\put(87,13){$x_{16}$}
\put(84,8){$T_{2222}$}

\end{overpic}
\end{center}
\caption{$2 \times 2 \times 2 \times 2$ tables}\label{table}
\end{figure}

For $2$ nowhere saturated $8$ dimensional faces, the set of extreme rays
of a cone is the set of $8$ columns of $16$ generators for the
semigroup, namely columns $1,2,7,8,11,12,13,14$.  
Since this $8$ dimensional face is spanned by $8$ vectors, it is a simplicial
face.  The extreme rays of the other cone is just the complement of
the extreme rays of this cone.
These two nowhere saturated faces are corresponding to $[234]$ marginals.
In order to see a picture, we
let $x_i$ correspond to the $i$th column of the defining
matrix $A$ above.  Let $T=(T_{ijkl})$ be a $2 \times 2 \times 2 \times 2$ table
where $i, j, k, l = 1, 2$.  Then we can write $16$ cells
of a $2 \times 2 \times 2 \times 2$ table as Figure \ref{table}.
There are eight $[234]$ marginals, i.e.
\begin{eqnarray*}
x_1+x_2 = T_{+111}, & x_3+x_4 = T_{+121},\\
x_5+x_6 = T_{+211}, & x_7+x_8 = T_{+221},\\
x_9+x_{10} = T_{+112}, & x_{11}+x_{12} = T_{+122},\\
x_{13}+x_{14} = T_{+212}, & x_{15}+x_{16} = T_{+222}.
\end{eqnarray*}
Four marginals $\{T_{+111}, T_{+221}, T_{+122}, T_{+212}\}$ correspond
to one face and and the other four marginals correspond to another
face.

We now consider $64$ almost saturated $2$ dimensional faces.
Each 2 dimensional face is spanned by two extreme rays and each extreme
ray corresponds to a cell $T_{ijkl}$.  Therefore each face can be
identified with a pair of cells.  
% we consider $[12][13][14][234]$ $2 \times 2 \times 2 \times 2$ model
% with cell denoted as     
Rather than listing all 64 faces, it is more instructive to consider
symmetry of the problem and list only different types of the faces.
Note that the semigroup has the symmetry with respect to interchanging
the values $1 \leftrightarrow 2$ for each $i,j,k,l$ independently, and
with respect to the permutation of indices
$j,k,l$. Therefore the product group $S_2 \times S_2 \times S_2 \times S_2
\times  S_3$ is naturally acting on the semigroup.  
By this action an almost saturated face is mapped to another almost
saturated face.
Detailed investigation of group invariance for Markov bases is given by
\cite{aoki-takemura-metr03-25} and \cite{aoki-takemura-metr05-14}.
The orbits for $2$ dimensional faces can be summarized as follows.
The first index $i$ can not be
interchanged with other $j$ or $k$ or $l$.
So we first look at $i$.  Then there are seven possible cases:
\begin{enumerate}
\item Whether we have a common value of $i$  or different
   values of $i$ in the two rays.
\item For each case, we can look at the number of common
   values of $j, \, k, \, l$ in the two rays.
\end{enumerate}
So there are seven possible cases as following:
\begin{enumerate}
\item $i$ is common and
\begin{enumerate}
\item\label{case1}  $j, \, k, \, l$ are all different.
\item \label{case2} one of $j, \, k, \, l$ is common.
\item \label{case3} two of $j, \, k, \, l$ are common.
\end{enumerate}
\item $i$ is different and
\begin{enumerate}
\item\label{case4}  $j, \, k, \, l$ are all different.
\item\label{case5}  one of $j, \, k, \, l$ is common.
\item\label{case6}  two of $j, \, k, \, l$ are common.
\item\label{case7}  $j, \, k, \, l$ are all common.
\end{enumerate}
\end{enumerate}

We counted the number of almost saturated $2$ dimensional faces
for each type among the seven cases and the results are the following:
\begin{itemize}
\item There are $8$ Type \ref{case1} almost saturated $2$ dimensional faces.
%\item There are $0$ Type \ref{case2} almost saturated $2$ dimensional faces. 
\item There are $24$ Type \ref{case3} almost saturated $2$ dimensional faces.
\item There are $8$ Type \ref{case4} almost saturated $2$ dimensional faces.
%\item There are $0$ Type \ref{case5} almost saturated $2$ dimensional faces. 
\item There are $24$ Type \ref{case6} almost saturated $2$ dimensional faces.
%\item There are $0$ Type \ref{case7} almost saturated $2$ dimensional faces.
\end{itemize}
There are no almost saturated $2$ dimensional faces of other types.
In order to make this classification based on symmetry clear, we give
the following example.

\begin{ex}
The cone generated by the $10$th and the $13$th columns of the matrix
and the cone generated by the $6$th and the $13$th columns of the matrix
are two dimensional almost saturated faces of the cone among $64$ cones.
Since the $6$th column represents $T_{2121}$, the $10$th column represents 
$T_{2112}$, and the $13$th column represents $T_{1122}$,
we have the index set for the cones:
\[
\begin{array}{|c|c|c|}\hline
&\mbox{1st cone} & \mbox{2nd cone}\\\hline
\mbox{1st ray}& 1122 & 1122\\
\mbox{2nd ray}& 2112 & 2121\\\hline
\end{array} 
\]
Note that for both cones, the value of $i$ in each ray has different value of
the other.  Thus we have Case 2.  Now we look at other indices.
$j$ is 1 in both cones.  However, $k$ has different values 1
and 2 in the first cone but $k$ has the same value in the second cone.  
$l$ has different values 1
and 2 in the second cone but $l$ has the same value in the first cone.
Thus the number of indices having the same value is two.  Therefore,  
the two cones are the same type, namely Type \ref{case6}.
In fact just by exchanging $k \leftrightarrow l$, the first cone is
mapped to the second cone in this example.
\end{ex}

\section*{Acknowledgment }

The authors would like to thank the referees for giving us helpful comments
to improve this paper.

%\newpage
\pagestyle{plain}
\bibliographystyle{plainnat}
\bibliography{faces}

\end{document}